
\magnification1200
\input amstex.tex
\documentstyle{amsppt}
\nopagenumbers
\hsize=12.5cm
\vsize=18cm
\hoffset=1cm
\voffset=2cm

\footline={\hss{\vbox to 2cm{\vfil\hbox{\rm\folio}}}\hss}

\def\DJ{\leavevmode\setbox0=\hbox{D}\kern0pt\rlap
{\kern.04em\raise.188\ht0\hbox{-}}D}

\def\txt#1{{\textstyle{#1}}}
\baselineskip=13pt
\def\hf{{\textstyle{1\over2}}}
\def\a{\alpha}\def\b{\beta}
\def\d{{\,\roman d}}
\def\e{\varepsilon}\def\E{{\roman e}}

\def\b{\beta} \def\g{\gamma}
\def\G{\Gamma}

\def\s{\sigma}

\def\={\;=\;}

\def\zt{\zeta(\hf+it)}

\def\D{\Delta}
\def\E{{\roman e}}

\def\r{\rho}
 
\def\z{\zeta}

\def\hf{{\textstyle{1\over2}}}
\def\txt#1{{\textstyle{#1}}}

\def\le{\leqslant} \def\ge{\geqslant}
\font\tenmsb=msbm10
\font\sevenmsb=msbm7
\font\fivemsb=msbm5
\newfam\msbfam
\textfont\msbfam=\tenmsb
\scriptfont\msbfam=\sevenmsb
\scriptscriptfont\msbfam=\fivemsb
\def\Bbb#1{{\fam\msbfam #1}}

\def \NN {\Bbb N}

\font\ff=cmr8
\def\txt#1{{\textstyle{#1}}}
\baselineskip=13pt

\font\teneufm=eufm10
\font\seveneufm=eufm7
\font\fiveeufm=eufm5
\newfam\eufmfam
\textfont\eufmfam=\teneufm
\scriptfont\eufmfam=\seveneufm
\scriptscriptfont\eufmfam=\fiveeufm
\def\mathfrak#1{{\fam\eufmfam\relax#1}}

\font\tenmsb=msbm10
\font\sevenmsb=msbm7
\font\fivemsb=msbm5
\newfam\msbfam
     \textfont\msbfam=\tenmsb
      \scriptfont\msbfam=\sevenmsb
      \scriptscriptfont\msbfam=\fivemsb
\def\Bbb#1{{\fam\msbfam #1}}

\def \NN {\Bbb N}

  \def\rightheadline{{\hfil{\ff
 Multiplicities of zeta zeors
 }\hfil\tenrm\folio}}

  \def\leftheadline{{\tenrm\folio\hfil{\ff
   Aleksandar Ivi\'c }\hfil}}
  \def\emptyheadline{\hfil}
  \headline{\ifnum\pageno=1 \emptyheadline\else
  \ifodd\pageno \rightheadline \else \leftheadline\fi\fi}

\topmatter

\title
On the  multiplicites of zeros of $\z(s)$ and its values over short intervals
\endtitle
\author   Aleksandar Ivi\'c
 \endauthor

\nopagenumbers

\medskip

\address
Aleksandar Ivi\'c, Katedra Matematike RGF-a
Universiteta u Beogradu, \DJ u\v sina 7, 11000 Beograd, Serbia
\endaddress

\bigskip
\email {
\tt
aleksandar.ivic\@rgf.bg.ac.rs, aivic\_2000\@yahoo.com }\endemail
\dedicatory
\enddedicatory
\abstract
{We investigate bounds for the multiplicities $m(\b+i\g)$, where
$\b+i\g\,$ ($\b\ge \hf, \g>0)$  denotes complex zeros of $\z(s)$. It is seen that the
problem can be reduced to the estimation of the integrals of the zeta-function
over ``very short'' intervals. A new, explicit  bound for $m(\b+i\g)$ is also derived, which
is relevant when $\b$ is close to unity. The related Karatsuba conjectures are also discussed.
 }
\endabstract
\endtopmatter

\document

\head
1. Introduction
\endhead

Let $r = m(\rho)$ denote the multiplicity of the complex
zero $\rho = \b + i\g$
of the Riemann zeta-function $\z(s)$.
It is defined for $\Re s >1$ by
$$
\z(s) \= \sum_{n=1}^\infty n^{-s},
$$
and otherwise by analytic continuation. This means that for some $r\in\NN$
$$
\z(\r) = \z'(\r)
= \ldots = \z^{(r-1)}(\r) = 0, \;\roman{but}\; \z^{(r)}(\r) \not = 0.
$$
All known zeros
$\r$ are simple (i.e., $m(\r) = 1$), and it may well be that they are all
simple, although the proof of this is
certainly beyond reach at present. Besides this strongest possible conjecture,
A.A. Karatsuba \cite{17} mentions two somewhat weaker conjectures: $m(\r) \ll 1\;(\forall\r)$
and $m(\r)$ is unbounded as $\g\to\infty$. He also says that the universality
of $\z(s)$ (see  S.M. Voronin \cite{26}) should include the last conjecture, but that all these
``are merely surmises''.

\medskip
In estimating $m(\r)$ one may suppose
that $\hf \le \b <1$ and that $\g > 0$,  since $\z(s)$ does not vanish
for $\Re s \ge1$, and $1-\r$ and $\overline\r$
 are zeros of $\z(s)$ if $\r$ is a zero. This follows from
 $\overline{\z(s)} = \z(\bar{s})$ and the functional equation
 $$
\z(s) \;=\; \chi(s)\z(1-s),\quad
\chi(s) := \frac{\G(\hf(1-s))}{\G(\hf s)}\pi^{s-1/2},
$$
where $\G(s)$ is the familiar gamma-function.
 For a comprehensive account on $\z(s)$,
 the reader is referred to the monographs of E.C. Titchmarsh \cite{23} and
 the author \cite{13}.

 \medskip
 Several results on the multiplicities of the zeros of the zeta-function were
 obtained in the author's paper \cite{14}. In particular, at the end of the paper
 it was stated that ``there is a possibility to bound $m(\b+i\g)$, provided one has a good
lower bound of the form
$$
\int_\delta^{2\delta}|\z(\b+i\g+i\alpha)|^k\d \alpha \,\ge\, \ell\,
= \,\ell(\g,\delta,k)\qquad(0 < \delta < {\txt{1\over4}},
\;\b \ge \hf,\;\g \ge \g_0 > 0) \leqno(1.1)
$$
for $k = 1,2$.'' Thus the problem is reduced to the evaluation of the moments of $\z(s)$
over ``very short'' intervals, namely integrals of the form
$$
\int_\delta^{2\delta}|\z(\b+i\g+i\alpha)|^k\d \alpha\qquad(0 < \delta < {\txt{1\over4}}),
\leqno(1.2)
$$
where $k\in\NN$ is fixed. The interval of integration can be justly
called ``very short'', since one assumes that $0 < \delta < {\txt{1\over4}}$.
One of the aims of this paper is to pursue further this approach and
analyze its potential.

\medskip
We note that zeta zeros with large multiplicities, statistically speaking,
are rare. Namely A. Fujii \cite{8} proved in 1975 that
$$
N_j(T) \;\le \;C_1N(T)\E^{-C_2\sqrt{j}}\qquad(T\ge T_0 >0),\leqno(1.3)
$$
where $N(T)$, as usual, denotes the number of complex zeros $\rho$ of $\z(s)$ for which
$0 < \Im \rho \le T$ (multiplicities counted), while $N_j(T)$
denotes those zeros counted by $N(T)$ whose
multiplicities are $j$. Here $j\,(\ge 1)$ is not necessarily fixed, and $C_1, C_2$ are
positive constants. A. Fujii \cite{9} in 1981 improved the exponential in (1.3) to $\exp(-C_2j)$, while
M.A. Korolev \cite{19} obtained much later in 2006 explicit numerical
values for the constants $C_1, C_2$
for the latter bound. Note that we have the identity
$$
\sum_{j=1}^\infty N_j(T) \;=\;N(T).\leqno(1.4)
$$
If $j = m(\b+i\g)$ with $0<\g\le T$, then by (2.11) one has $j \ll \log\g$,
hence it follows that the sum in (1.5) is finite.

It seems plausible that uniformly, for any given $j\ge2$,
$$
N_j(T) \;=\; o\bigl(N(T)\bigr)\qquad(T\to\infty), \leqno(1.5)
$$
which implies that $N(T)\sim N_1(T)$, namely that all the zeros are simple.
However, in general, (1.5) is not known yet. It follows from (1.3) if $j = j(T)\to\infty$ as $T\to\infty$.
 The bound in (1.3) suggests that $N_j(T)$ is
a non-increasing function of $j$ for a fixed $T$, but this is not easy to prove.
Note that the relation (1.5) certainly cannot hold for $j=1$, since D.R. Heath-Brown \cite{12}
showed that $N_1(T) \gg N(T)$.

\medskip
In Section 2 and Section 4 we shall deal with lower bounds of the form (1.1) and obtain
in Theorem 2 a new lower bound. In Section 3 we shall consider the Karatsuba conjectures
involving the quantity
$$
F(T,\D) := \max_{t\in[T,\, T+\D]} |\zt|\qquad(0 < \D \le1),
$$
which is closely related to the integral in (1.1). Finally,
 in Section 5 we shall employ a complex integration technique to
obtain an explicit upper bound for $m(\b+i\g)$, which is relevant when $\b$ is close to unity.
\medskip
\head
2. Integrals over short intervals
\endhead
\medskip
The argument for the estimation of $m(\rho)=r$ that leads to (1.1) is as follows.
 For fixed $\b$ such that $\b \ge \hf$, let $\Cal D$ be the
rectangle with vertices
$$
\textstyle {1\over4} -\b\pm i\log^2\g,\; 2 \pm i\log^2\g,\; \z(\rho) = 0,\;
\rho = \b+i\g \;\;(\g \ge \g_0 > 0),
$$
and let $\a$ be a parameter for which $0 < \alpha \le 1$. Since $\rho$ is a zero of $\z(s)$
of multiplicity $r$, the function $\z(s+\rho)s^{-r}$ is regular at $s=0$.
By the residue  theorem we obtain
$$
{\z(\b+i\g+i\alpha)\over(i\alpha)^r}
= {1\over2\pi i}\int_{\Cal D} \Gamma(s-i\alpha){\z(s+\rho)\over s^r}\d s.\leqno(2.1)
$$
Namely of the poles of the gamma-factor only $s = i\alpha\,$ is  in $\Cal D$, and it is a simple pole.
The unique pole of $\z(s+\r)$, namely $s = 1-\r$, lies outside $\Cal D$.
This gives, in view of the fast decay of the gamma-function (see e.g., (A.34) of \cite{13}),
$$
\z(\b+i\g+i\alpha) \ll \alpha^r\left(\g(\b - {\textstyle {1\over4}})^{-r}
+ 2^{-r}\right) \ll \alpha^r\g{(\b - {\textstyle {1\over4}})}^{-r},\leqno(2.2)
$$
and the case when
$\g{(\b - {\textstyle {1\over4}})}^{-r} \ll 2^{-r}$ is easy, since it implies that
$$
\g \ll (\b - {\textstyle {1\over4}})^{r}2^{-r} \le \left(\frac38\right)^r,
$$
and this is impossible if $r \ge r_0$. Hence either $2^{-r} \ll \g(\b - {\textstyle {1\over4}})^{-r}$
or
$$
r \;=\; m(\b + i\g) \;\ll\; 1,
$$
and this case is covered by the term $O(1)$ in (2.3). It is, of course,
possible to insert in the  integrand in (2.1) the factor $X^{s-i\a}\;(X>1)$, and try to use
convexity. This does not appear to give any substantial improvement.
Consequently, if $\delta$ is a constant satisfying $0 < \delta< {1\over8}$,
then raising (2.2) to the power $k$ and integrating over $\a$ we have
$$
\int_\delta^{2\delta}|\z(\b+i\g+i\alpha)|^k\d \alpha
\ll \g^k{(\b - {\textstyle {1\over4}})}^{-rk}\int_\delta^{2\delta}\alpha^{rk}
\d\alpha \ll (8\delta)^{rk}\g^k.
$$
Thus, recalling (1.1) and taking logarithms, we have

\medskip
{\bf Theorem 1}. {\it If $\b\ge\hf,\, \g>\g_0>0,\,k>0\,$ and $ 0 < \delta < \txt{1\over8}\,$,
 then with the notation introduced above we have}
$$
m(\b+i\g) = r \le {1\over\log
\left({1\over8\delta}\right)}
\left(\log\g - {1\over k}\log\ell + O(1)\right) + O(1).\leqno(2.3)
$$

\medskip
Therefore (2.3) shows that the upper bound for $m(\b+i\g)$ can be made to
depend on $\ell$ in (1.1), that is, on lower bounds for moments of $\z(s)$ over very
short intervals. We would like to let $\delta \to 0+$ in (2.3)
and obtain
$$
m(\b+i\g) \;=\;o(\log \g)\qquad(\b\ge\hf,\;\g\to\infty).\leqno(2.4)
$$
This relation is equivalent to
$$
\lim_{\delta\to0+}\,\frac{\ell}{\log\Bigl({1\over8\delta}\Bigr)} =
\lim_{\delta\to0+}\,\frac{\ell(\g,\delta,k)}{\log\Bigl({1\over8\delta}\Bigr)}\;=0.
$$

 However, by using the argument on top of p. 219 of E.C. Titchmarsh \cite{23}
and the first  inequality on p. 230, it follows that
(1.1) holds with $\ell = \delta \gamma^{-A/\delta}$.
By suitably elaborating the method it follows that even
$$
\ell \;=\; \delta \gamma^{A\log\delta}\leqno(2.5)
$$
is permissible, for some absolute $A > 0$. These bounds, unfortunately, are
too weak to yield (2.4). The bound in (2.5) can be compared to  the case $\s = \hf$ of
Theorem 2 in Section 4.

\medskip
We remark that, on the Lindel\"of Hypothesis (LH) that $\z({1\over2} +it) \ll_\e
 |t|^\varepsilon$, one has indeed (2.4). Note that $f(x) \ll_{\a,\b,\ldots} g(x)$ (same
 as  $f(x) = O_{\a,\b,\ldots}\left\{g(x)\right\}$) means that the implied $\ll$--constant
 (resp. $O$-constant) depends on $\a,\b,\ldots\,$.
Also on the Riemann Hypothesis (RH, well-known that it implies the LH; see \cite{23})
 that $\r = {1\over2} + i\g\;(\forall \rho)$ one has
$$
m(\b + i\g) \;\ll\; {\log \g\over\log\log\g}.\leqno(2.6)
$$
Furthermore, on the RH, H.L. Montgomery \cite{21} proved
that at least 2/3 of the zeta zeros are simple, namely $N_1(T) \ge \frac23 N(T)\;(T\ge T_0)$. 
His result was recently improved by
H.M. Bui and D.R.Heath-Brown \cite{2} (also on the RH), who obtained
the constant $19/27 = 0.\overline{703}\,$ in place of 2/3.

\medskip
It transpires that the estimation of $m(\b + i\g)$ is a very difficult
problem, and one which is not satisfactorily solved even under the
assumption of the LH or the RH.
To see how one obtains (2.4) and (2.6) recall that
for $N(T)$, the number of zeros $\b+ i\g$ for which $0 < \g \le T$, one
has the classical Riemann-von Mangoldt formula (see \cite{13} or \cite{23} for a proof)
$$
N(T) = {T\over2\pi}\log\left({T\over2\pi}\right) - {T\over2\pi} + {7\over8}
+ S(T) + O\left({1\over T}\right),
\leqno(2.7)
$$
where $S(T) = {1\over\pi}\arg\z(\hf+ iT),$ and the term $O(1/T)$ is a smooth function.
Here $\arg\z({1\over2}+ iT)$ is obtained by continuous
variation along the segments
joining the points $2, 2 + iT, {1\over2} + iT$, starting with the value 0.
If $T$ is the ordinate of a zero lying on the critical line,
 then $S(T) =S(T + 0)$. One has (see \cite{23}) the bounds
$$
S(T) \ll \log T,\quad S(T) = o(\log T) \quad({\roman {LH}}),\quad
S(T) \ll {\log T\over\log\log T}\quad ({\roman {RH}}).
$$
These bounds combined with (2.7) and the trivial inequality
$$
m(\b + i \g) \;\le\; N(\g + H) - N(\g - H)  \qquad(0 < H \le 1)\leqno(2.8)
$$
easily yield
$$
m(\b + i \g) \ll \log\g,\quad m(\b + i \g) = o(\log\g) \quad({\roman {LH}}),\quad
m(\b + i \g) \ll \frac{\log\g}{\log\log\g}\quad ({\roman {RH}}),
$$
respectively.
It seems, however, that these estimates
are much too large, and that perhaps one even has
$$
m(\b + i\g) \;\ll_\varepsilon \;(\log\log\g)^{1+\varepsilon},\leqno(2.9)
$$
which is  weaker than the previously stated conjectures, in particular that all zeros
are simple.
The direct use of pointwise estimates for $S(T)$ certainly cannot give anything
close to (2.9), since one has
$$
S(T) = \Omega_\pm\left(\left({\log T\over\log\log T}\right)^{1/3}\right),
\quad
S(T) = \Omega_\pm\left(\left({\log T\over\log\log T}\right)^{1/2}\right)\quad
({\roman {RH}}).\leqno(2.10)
$$
This was proved proved by K.-M. Tsang \cite{25} (his result is unconditional)
and H.L. Montgomery \cite{22}, respectively.
As usual, $f(x) = \Omega_\pm\bigl(g(x)\bigr)$ means that the inequalities
$$
\limsup_{x\to\infty}\frac{f(x)}{g(x)}
>0 \quad\roman{and}\quad \liminf_{x\to\infty}\frac{f(x)}{g(x)} <0
$$
both hold. One could use (2.8) with $H = o(1)\;(\g\to\infty)$ to
try to improve the existing bound
$$
m(\b+i\g) \;\ll\;\log\g\qquad(\hf \le \b < 1).\leqno(2.11)
$$
 In view of (2.7) this is equivalent to
obtaining bounds for $S(\g + H) - S(\g - H)$, but no satisfactory
results seem to be known for this problem. Note that (2.11) easily
follows from (2.7), (2.8) and $S(T) \ll \log T$.
In spite of all the efforts, this is still the
best unconditional bound for the whole range $\hf \le \b < 1$.
For an additional discussion concerning $S(T)$, see Section 6.

\medskip
\head
3. The Karatsuba conjectures
\endhead
\medskip
A function closely related to the integral in (1.1) (when $\b = \hf,\,k=1$) is
$$
F(T,\D) := \max_{t\in[T,\, T+\D]} |\zt|\qquad(0 < \D \le1),\leqno(3.1)
$$
where $\D$ may depend on $T$.
Namely, for a fixed $k>0$, one clearly has
$$
\int_\delta^{2\delta}|\z(\hf+i\g+i\alpha)|^k\d \alpha
\;=\; \int_0^{\delta}|\z(\hf+i\g+i\delta+ix)|^k\d x
\;\le\; \delta F^k(\g +\delta, \delta).
\leqno(3.2)
$$
The quantity $F(T,\D)$ was introduced and studied by A.A. Karatsuba \cite{15},
\cite{16}, \cite{17}. He made the following conjectures.

\medskip
{\bf Conjecture 1}. There exists a positive function $\D = \D(T) \to0$ as $T\to\infty$ such that,
for some constant $A>0$,
$$
F(T,\D) \;\ge \; T^{-A}.\leqno(3.3)
$$
\medskip
{\bf Conjecture 2}.
Conjecture 1 is valid for $\D = (\log\log T)^{-1}$.

\medskip
{\bf Conjecture 3}.
Conjecture 1 is valid for $\D = (\log T)^{-1}$.

\medskip
These conjectures have not been proved unconditionally yet.
Clearly Conjecture 3  implies
Conjecture 2, which in turn implies Conjecture 1. M. Garaev \cite{10} proved that the RH
implies Conjecture 3, while Karatsuba himself showed unconditionally that
$$
F(T,\D) \;\ge\; \E^{A\log\D\log T}\qquad(0 < \D \le 1/(\log T)).\leqno(3.4)
$$
Shao-Ji Feng \cite{6} proved that the LH implies Conjecture 1 with an arbitrary constant
$A>0$. Other relevant works on this subject include the papers of M.E. Changa \cite{5},
B. Kerr \cite{18} and M.A. Korolev \cite{20}.

\medskip
In view of (3.1) and (3.2) it is seen that the Karatsuba conjectures have their counterparts
involving the integral in (3.2). For example, the conjecture
$$
\int_\delta^{2\delta}|\z(\hf+iT+i\alpha)|^k\d \alpha \;\gg\; T^{-A}\qquad(\delta = \delta(T)\to 0)
$$
is less stringent than Karatsuba's Conjecture 1, and similarly for the other two conjectures.

\medskip
We have

\medskip
{\bf Theorem 2}. {\it If Conjecture 1 holds, then}
$$
m(\hf + i\g) \;=\;o(\log \g)\qquad(\g\to\infty). \leqno(3.5)
$$

\medskip
{\bf Proof}. The assertion follows from (2.2) with $\b = \hf, \a = \D = \D(\g)$. Namely
Conjecture 1 gives
$$
\g^{-A} \;\ll\; \D^r\g{(\b - \txt{\frac{1}{4}})}^{-r} \;=\; ({4\D})^r\g,
$$
which implies
$$
\left(\frac{1}{4\D}\right)^r \;\ll\;\g^{A+1}.
$$
Taking logarithms, we obtain
$$
m(\hf+i\g)\log\frac{1}{4\D} = r \log\frac{1}{4\D} \le C + (A+1)\log \g\qquad(\g \ge \g_0 >0),
$$
and the assertion readily follows, since
$$
\lim_{\g\to\infty}\log\,\frac{1}{4\D} \;\to\;+\infty
$$
by the assumption on $\D = \D(\g)$ in Conjecture 1. This shows again that the LH implies
$m(\hf+i\g) = o(\log\g)$. A conditional result, similar to Theorem 2, is given by A.A. Karatsuba
\cite{17}. Naturally, Conjecture 2 and Conjecture 3, with explicit values
of $\D = \D(T)$ would lead to sharper results on $m(\hf + i\g)$. Open questions are:
does (3.5) imply the LH or Conjecture 1?

\medskip
\head
4. Integrals of $|\z(\s+it)|$ over very short intervals
\endhead

\medskip
We have the following result, which is more general than Karatsuba's bound (3.4),
but of the same strength. The method of proof is different from Karatsuba's.

\medskip
{\bf Theorem 3}. {\it For $k>0, \hf \le \s \le 1, 0 < \delta \le \hf, T\ge T_0 >0$
and a suitable constant $C>0$ we have}
$$
\int_{T-\delta}^{T+\delta}|\z(\s+it)|^k\d t \;\ge\; 2\delta T^{-Ck\log(\E/\delta)}.
\leqno(4.1)
$$

\medskip
{\bf Proof}. We start from Th. 9.6 (B)  of Titchmarsh's
book \cite{23}, namely from the classical formula
$$
\log \z(s) = \sum_{|t-\g|\le1}\log(s-\rho) + O(\log t),
$$
which is valid unconditionally
 for $-1\le\s\le2, s\ne \rho, -\pi < \Im \log(s-\rho) \le \pi$, where $\rho$
 denotes complex zeros of $\z(s)$.
Since $\Re \log z = \log|z|$, then by taking real parts in this formula it follows that
$$\eqalign{
\log|\z(s)| & = \sum_{|t-\g|\le1}\log|s-\rho| +O(\log t)\cr&
\ge \sum_{|t-\g|\le1}\log|t-\g| + O(\log t).\cr}\leqno(4.2)
$$
To get rid of the logarithms one uses (this is a consequence of the arithmetic-geometric
means inequality)
$$
\log\Bigl\{{1\over b-a}\int_a^b f(t)\d t\Bigr\}
 \ge {1\over b-a}\int_a^b\log f(t)\d t\leqno(4.3)
$$
for $a < b, f(t) \in L[a,b]$ and $f(t)>0$ in $[a,b]$. Hence with
$$
a = T-\delta,\; b = T+\delta,\; f(t) =|\z(\s+it)|^k,
$$
(4.3) yields
$$
\log\Bigl\{{1\over2\delta}\int_{T-\delta}^{T+\delta}
|\z(\s+i\a)|^k\d \a\Bigr\} \ge{1\over2\delta}\int_{T-\delta}^{T+\delta}k\log|\z(\s+it)|\d t.
\leqno(4.4)
$$
Note that we have
$$
\eqalign{&
\int_{T-\delta}^{T+\delta}\sum_{|t-\g|\le1}\log|t-\g|\d t\cr&
= \int_{T-\delta}^{T+\delta}\sum_{|t-\g|\le\delta}\log|t-\g|\d t
+ \int_{T-\delta}^{T+\delta}\sum_{\delta<|t-\g|\le1}\log|t-\g|\d t\cr&
= I_1 + I_2,\cr}
$$
say. But, since $\log|t-\g|\le0$ for $|t-\g|\le1$ and
$$
[\max(T-\delta,\g-\delta),\,\min(T+\delta,\,\g+\delta)] \;\subseteq\;
[\g-\delta, \g+\delta],
$$
we obtain
$$
\eqalign{
I_1& = \sum_{T-2\delta\le\g\le T+2\delta}\int_{\max(T-\delta,\g-\delta)}^
{\min(T+\delta,\g+\delta)}
\log|t-\g|\d t\cr&
\ge \sum_{T-2\delta\le\g\le T+2\delta}\int_{\g-\delta}^{\g+\delta}\log|t-\g|\d t\cr&
= \sum_{T-2\delta\le\g\le T+2\delta}\int_{-\delta}^\delta \log|u|\d u\cr&
= 2(\delta\log\delta -\delta)\sum_{T-2\delta\le\g\le T+2\delta}1\cr&
\ge -C\delta\log({\E}/\delta)\log T,\cr}
$$
since $\delta\log\delta -\delta<0$ for $0<\delta\le 1$. We also have,
since $S(T) \ll \log T$ and $\log\delta<0$,
$$
I_2 = \int\limits_{T-\delta}^{T+\delta}\sum_{\delta<|t-\g|\le1}\log|t-\g|\d t
\ge \log\delta\int\limits_{T-\delta}^{T+\delta}\sum_{\delta<|t-\g|\le1}1\d t
\ge -C\delta\log(\E/\delta)\log T,
$$
where $C$ is a positive constant.
Therefore from (4.2), (4.4) and the above bounds we obtain
$$
\log\Bigl\{{1\over2\delta}\int_{T-\delta}^{T+\delta}
|\z(\s+i\a)|^k\d \a\Bigr\} \ge -kC\log(\E/\delta)\log T,
$$
which implies the lower bound in Theorem 3. This completes the proof.
We remark that (4.1) in conjunction with (2.3) produces only the classical
bound (2.11).

\medskip
{\bf Remark 1}. Note that Karatsuba's function $F(T,\D)$ (see (3.1)) can be connected to
the integral of $\log|\zt|$ over a very short interval. Namely, for $0 < \D \le 1$,
using (4.3) we have
$$
\eqalign{&
F(T,\D) = \max_{0\le u\le\D}|\z(\hf + iT + iu)|\cr&
\ge \frac1\D\int_0^\D |\z(\hf + iT + iu)|\d u \ge \exp\left\{\frac1\D \int_0^\D
\log|\z(\hf + iT + iu)|\d u\right\}.\cr}
$$
Putting $T_0 = T + \hf\D, \delta = \hf\D$, it follows that
$$
F(T,\D) \ge \exp\left\{\frac{1}{2\delta} \int_{T_0-\delta}^{T_0+\delta} \log|\zt|\d t\right\}
\quad(0 < \delta \le \hf).\leqno(4.5)
$$
The integral in (4.5) is precisely of the type that was dealt with in the proof of Theorem 3.

\medskip
{\bf Remark 2}. Note that if (4.1) is known to hold for $k=1$, then one can easily
deduce that it holds for $k>1$ as well. Namely, by H\"older's inequality for
integrals we have, for $k>1$,
$$
\int_{T-\delta}^{T+\delta}|\z(\s+it)|\d t \le
\left(\int_{T-\delta}^{T+\delta}|\z(\s+it)|^k\d t\right)^{1/k}(2\delta)^{1-1/k}.
\leqno(4.6)
$$
Therefore if
$$
\int_{T-\delta}^{T+\delta}|\z(\s+it)|\d t \;\ge\; 2\delta T^{-C\log(\E/\delta)},
$$
one easily obtains (4.1) from (4.6).

\bigskip
\head
5. A bound for multiplicities when $\b$ is close to unity
\endhead
\medskip
We finally present an explicit bound for $m(\b+i\g)$, which is relevant when $\b$ is close to unity.
If such $\b$ exists, then the RH cannot hold.
The result  is

\bigskip
{\bf Theorem 4}. {\it Let $5/6 \le \b < 1$. Then we have, for $\g\ge\g_0(\e)$,
a suitable constant $C>0$ and any $\e>0$,}
$$
\eqalign{
m(\b+i\g) &\le C + \frac{13.35\b}{3(1-\b)\log 6 + \b\log 2}(1-\b)^{3/2}\log \g\cr&
+ \frac{7(3-2\b)+\e}{9(1-\b)\log6+3\b\log2}\log\log\g.\cr}\leqno(5.1)
$$
\bigskip
{\bf Corollary 1}.   For $5/6 \le \b < 1$ and $\g\ge\g_1>0$, we have
$$
m(\b+i\g) \;\le\;4\log\log\g + 20(1-\b)^{3/2}\log \g.\leqno(5.2)
$$
\medskip
{\bf Corollary 2}. If $m(\b+i\g) \ge 8\log\log\g$ for $5/6 \le \b < 1$ and $\g\ge\g_2>0$, then
$$
\b  \;\le\; 1 - {\left(\frac{m(\b+i\g)}{40\log\g}\right)}^{2/3}.\leqno(5.3)
$$

\medskip
One obtains (5.2) and (5.3) by noting that
$$
\frac{13.35\b}{3(1-\b)\log 6 + \b\log 2} \le \frac{13.35}{\log2} = 19.25997\ldots\,
$$
and that $m(\b+i\g) \ge 8\log\log\g$ implies  $4\log\log\g \le \hf m(\b+i\g)$. The
bound (5.3) says that, if the zero $\b+i\g$ has a large multiplicity, then $\b$
cannot be large.

\medskip
{\bf Proof of Theorem 4}. This result is a sharpening of Theorem 4 of \cite{14}, where
one had the Vinogradov symbol $\ll$ instead of explicit inequalities.
Let $\b \ge 5/6$, $r = m(\b+i\g)$ and $\Cal E$ be the rectangle
with vertices $-2(1-\b) \pm 2i\log^2\g,\,1\pm 2i\log^2\g$.
If $X\,(0 < X \ll \g^C)$ is a parameter which will be suitably chosen,
then by the residue theorem we obtain
$$
{\z(1-\b+\r)\over(1-\b)^r} \;=\; {1\over2\pi i}
\int_{\Cal E}X^{s-1+\b}\G(s-1+\b){\z(s+\r)\over s^r}\d s
\quad(\rho = \b+i\g),\leqno(5.4)
$$
which is similar to (2.1).
Namely $-\b < -2(1-\b) < 1-\b$, while $\G(s-1+\b)$ has simple
poles at $s = 1-\b, -\b, -1-\b,\ldots$.
For the gamma-function we shall use the estimate
$$
\G(w)  \;\ll\; {\E^{-|{\roman {Im}}\,w|}\over |w|}.
$$
To bound the zeta-factor on the left side of  (5.4) we shall use the inequality
$$
|\z(\s+it)| \;\le \; At^{B(1-\s)^{3/2}}\log^{2/3}t\qquad( t\ge 3,\quad \hf \le \s \le 1),
\leqno(5.5)
$$
with the currently best known values $A = 76.2, B = 4.45$, due to K. Ford \cite{7}.
For our purposes it is the value of the constant $B$ that is relevant. On the left side of
$\Cal E$ we have $\Re(s+\r) = 3\b - 2 \ge 1/2$, since $\b \ge 5/6$ is assumed to hold.

We shall also use the bound
$$
\z(1+it) \;\gg\;{(\log |t|)}^{-2/3}{(\log\log |t|)}^{-1/3},
$$
which is a consequence of Lemma 12.3 of \cite{13}. Like (5.5), this bound is obtained by an elaboration
of the classical method of Vinogardov--Korobov (see e.g., Chapter 6 of \cite{13}) for the estimation
of certain exponential sums.
It follows then from (5.4) that
$$
\eqalign{&
{(1-\b)^{-r}\over\log^{2/3}\g(\log\log\g)^{1/3}} \ll \E^{-\log^2\g} + X^\b\cr&+ 2^{-r}(1-\b)^{-r}
X^{-3(1-\b)}\log\g\,\max_{|t|\le\log^2\g}\,|\z(3\b-2+i\g+it)|.\cr}\leqno(5.6)
$$
Using (5.5) in (5.6) it follows that
$$
(1-\b)^{-r} \ll L(\g)\left(X^\b + 2^{-r}(1-\b)^{-r}X^{-3(1-\b)}\g^{3B(1-\b)^{3/2}}\right),\leqno(5.7)
$$
where for brevity we put
$$
L(\g) := (\log\g)^{7/3}(\log\log\g)^{1/3}.
$$
We multiply (5.7) by $2^r(1-\b)^r$ and use $1-\b\le 1/6$ to deduce  that
$$
2^r \ll \left(3^{-r}X^\b + X^{-3(1-\b)}\g^{3B(1-\b)^{3/2}}\right)L(\g).\leqno(5.8)
$$
Now we choose $X$ in (5.8) so that the two terms on the right-hand side are equal. Thus
$$
X = 3^{r/(3-2\b)}\g^{3B(1-\b)^{3/2}/(3-2\b)} \qquad(\ll \g^C).
$$
This gives
$$
2^r \ll 3^{-r}3^{\b r/(3-2\b)}\g^{3B\b(1-\b)^{3/2}/(3-2\b)}L(\g).
$$
We raise this to the power $3-2\b$ and take logarithms to obtain
$$
\eqalign{&
r(3-2\b)\log 2 + r(3-2\b)\log 3 - \b r\log 3 \cr &
\le C_1 + 3B\b(1-\b)^{3/2}\log\g + (3-2\b)\log L(\g).\cr}\leqno(5.9)
$$
Since the coefficient of $r$ on the left-hand side equals
$$
3(1-\b)\log 6 + \b\log 2,
$$
and
$$
(3-2\b)\log L(\g) \le \frac{7(3-2\b)+\e}{3}\log\log\g,
$$
we obtain the assertion (5.1) of Theorem 4 from (5.9).

\bigskip
\head
6. Some remarks concerning $S(T)$
\endhead
\bigskip
We conclude with some remarks concerning the function $S(T)$ and its effects on the estimation
of $m(\b+i\g)$.
In the paper of Goldston--Gonek \cite{11} it is proved, under the RH, that
$$
|S(T+H) - S(T)| \le \left({1\over2} + o(1)\right){\log T\over\log\log T}
\quad(T\to\infty, 0 < H \le \sqrt{T}\,).\leqno(6.1)
$$
This implies, under the RH, in view of (2.7) and (2.8), the explicit upper bound
$$
m(\b + i\g) \;\le\; \left(\frac{1}{2} + o(1)\right)\frac{\log\g}{\log\log\g}
\qquad(\hf\le\b< 1,\;\g\to\infty),\leqno(6.2)
$$
on taking $H = 1/\log^2\g$, say.
It is known that, unconditionally (see E.C. Titchmarsh \cite{23}) one has,
$$
\int_0^T S(t)\d t \;\ll\; \log T.\leqno(6.3)
$$
From (6.3) it follows  that every interval $[T, T+\log^2T]$ contains a point
$t_0$ for which $S(t_0) \le 1$,
and a point $t_1$ for which $S(t_1) \ge -1$.
From this and (6.1) one obtains
$$
S(T) \le \left(\frac{1}{2} + o(1)\right)\frac{\log T}{\log\log T}\qquad(\roman{RH}, \;T\to\infty).
\leqno(6.4)
$$
The constant one half in (6.4) (and thus also in (6.2)) was improved by Carneiro, Chandee
and Milinovich \cite{3} to 1/4, and the ``o(1)" term
is actually $$O\Bigl(\frac{\log\log\log T}{\log\log T}\Bigr).$$
Generalizations of (6.4)
to suitable $L$-functions were recently established in
a paper by E. Carneiro and R. Finder \cite{4}.

\medskip
A recent unconditional, explicit bound for $S(T)$ is
$$
|S(T)| \le 0.111\log T + 0.275\log\log T + 2.450,
$$
which is valid for $T \ge \E$. This is a recent result of T. Trudgian \cite{24}. By (2.7) and (2.8)
it immediately implies the unconditional bound
$$
m(\b+i\g) \;\le   2(0.111\log \g + 0.275\log\log \g + 2.450)\qquad(\hf \le \b <1, \g \ge 14),
$$
which is an explicit version of (2.11).

\medskip
The largest known values of $S(T)$ (in absolute value) at present are are,
for $T$ less than 29 trillion ($\approx$ means approximately):
$$
S(T) \approx 3.0214, \; T \approx 53\, 365\, 784\, 979;\; S(T)
\approx -3.2281,\; T \approx 69\, 976\, 605\, 145.
$$
This was found by S. Wedeniwski \cite{27} and his team in the larger context of searching for the zeros
of $\z(s)$ on the critical line. The first 100 billion zeros are simple and lie on the critical line.
More extensive calculations are to be found in the forthcoming paper of J.W. Bober and G.A. Hiary [1].
This shows that the values of $T$ needed for the $\Omega$-results in (2.10) to take effect must be
extremely large.

\vfill
\eject
\topglue1cm
\bigskip
\Refs
\bigskip

\item{[1]} J.W. Bober and G.A. Hiary, New computations of the Riemann zeta function on the critical
line, to appear, preprint available at {\tt arXiv:1607.00709}.

\smallskip
\item{[2]}
H.M. Bui and D.R. Heath-Brown, On simple zeros of the Riemann zeta-function,
 Bull. Lond. Math. Soc. {\bf45}(2013), no. 5, 953-961.

\smallskip
\item{[3]} E. Carneiro,  V. Chandee and
M.B. Milinovich, Bounding $S(t)$ and $S_1(t)$ on the Riemann hypothesis,
Math. Ann. {\bf356}(2013), 939-968.

\smallskip
\item {[4]} E. Carneiro and R. Finder, On the argument of $L$-functions,
Bull. Braz. Math. Soc. (N.S.) {\bf46}(2015), no. 4, 601-620.

\smallskip
\item{[5]} M.E. Changa, Lower bounds for the Riemann zeta function on the critical line,
Math. Notes {\bf76}(2004), 859-864.

\smallskip
\item{[6]} Shao-Ji Feng, On Karatsuba conjecture and the Lindel\"of hypothesis,
Acta Arithmetica \break {\bf114}(2004), 295-300.

\smallskip
\item{[7]} K. Ford, Vinogradov's integral and bounds for the Riemann zeta function,
Proc. Lond. Math. Soc. (3){\bf85}(2002), 565-633.

\smallskip
\item{[8]} A. Fujii, On the distribution of the zeros of the Riemann zeta-function in short
intervals, Bull. Amer. Math. Soc. {\bf 81}(1975), 139-142.

\smallskip
\item{[9]} A. Fujii, On the zeros of Dirichlet $L$-functions. II.
(With corrections to ``On the zeros of Dirichlet $L$-functions. I'' and the subsequent papers),
Trans. Amer. Math. Soc. {\bf267}(1981), 33-40.

\smallskip
\item{[10]} M.Z. Garaev, Concerning the Karatsuba conjectures, Taiwanese
J. Math. {\bf6}(2002), 573-580.

\smallskip
\item{[11]}
D.A.  Goldston and S.M. Gonek,  A note on $S(t)$ and the zeros of the Riemann zeta-function,
 Bull. Lond. Math. Soc. {\bf39}(3)(2007),  482-486.

\smallskip
\item{[12]} D.R. Heath-Brown, Simple zeros of the Riemann zeta-function, Bull. London Math. Soc.
{\bf11}(1979), 17-18.

\smallskip
\item{[13]} A. Ivi\'c, The Riemann zeta-function, John Wiley \&
Sons, New York 1985 (reissue,  Dover, Mineola, New York, 2003).

\smallskip
\item{[14]} A. Ivi\'c,  On the multiplicity of  zeros of the zeta-function,
Bulletin CXVIII de l'Acad\'e\-mie Serbe des Sciences et des Arts - 1999,
Classe des Sciences math\'ematiques et naturelles,
Sciences math\'ematiques No. {\bf24}, pp. 119-131.

\smallskip
\item{[15]} A.A. Karatsuba, On lower estimates of the Riemann zeta-function,
Dokl. Akad. Nauk {\bf376}\break (2001), 15-16.

\smallskip
\item{[16]} A.A. Karatsuba, Lower bounds for the maximum modulus of $\z(s)$
in small domains of the critical strip, Math. Notes {\bf70}(2001), 724-726.

\smallskip
\item{[17]} A.A. Karatsuba, Zero multiplicity and lower bound estimates of $|\z(s)|$,
Funct. Approx. Comment. Math. {\bf35}(2006), 195-207.

\smallskip

\item{[18]} B. Kerr, Lower bounds for the Riemann zeta function on short
intervals of the critical line, Archiv Math. {\bf105}(2015), 45-53.

\smallskip
\item{[19]} M.A. Korolev, On multiple zeros of the Riemann zeta-function,
Izv. Math. {\bf70}(2006), 427-446; translation from
Izv. Ross. Akad. Nauk, Ser. Mat. {\bf70}(2006), 3-22.

\smallskip
\item{[20]} M.A. Korolev,
On large values of the Riemann zeta-function on short segments of the critical line,
Acta Arith. {\bf166}(2014), 349-390.

\smallskip
\item{[21]} H.L. Montgomery,  The pair correlation of zeros of the
zeta-function, Proc. Symp. Pure Math. {\bf 24}, AMS, Providence R.I.,
1973, 181-193.

\smallskip
\item{[22]} H.L. Montgomery,  Extreme values of the Riemann
zeta-function, Comment. Math. Helv. {\bf 52}(1977), 511-518.

\smallskip
\item{[23]}  E.C. Titchmarsh, The theory of the Riemann
zeta-function (2nd edition),  Oxford University Press, Oxford, 1986.

\smallskip
\item{[24]} T. Trudgian. An improved upper bound for the argument of the
Riemann zeta-function on the critical line II, J. Number Theory {\bf134}(2014), 280-292.

\smallskip
\item{[25]} K.-M. Tsang,  Some $\Omega$--theorems for the Riemann
zeta-function, Acta Arithmetica {\bf 46}(1986), 369-395.

\smallskip
\item{[26]} S.M. Voronin, A theorem on the ``universality'' of the Riemann zeta-function,
Izv. Akad. Nauk SSSR Ser. Mat. {\bf39}(1975), no. 3, 475-486.

\smallskip
\item{[27]}
S. Wedeniwski, Results connected with the first 100 billion zeros of the Riemann
zeta function, 2002, at {\tt http://piologie.net/math/zeta.result.100billion.zeros.html}

\endRefs
\vskip2cm

\enddocument

\bye